\documentclass[12pt]{article}
\usepackage{amsmath}
\usepackage{amssymb}
\usepackage{theorem}
\usepackage{a4}
\usepackage{pb-diagram}
\usepackage{lamsarrow}
\usepackage{pb-lams}

\def\gR{{\mathbb R}}

\def\gC{{\mathbb C}}

\def\cF{{\cal F}}
\def\cG{{\cal G}}
\def\cH{{\cal H}}

\def\cN{{\cal N}}

\def\cZ{{\cal Z}}

\def\tX{{\mathfrak X}}

\def\tU{{\mathfrak U}}
\def\tsl{{\mathfrak{sl}}}
\def\tsu{{\mathfrak{su}}}

\def\fA{{\mathfrak A}}
\def\fB{{\mathfrak B}}

\def\findem{~\hfill$\square$}

\def\demo{\noindent {\sl Proof: }}

\def\der{{\mathrm{Der}}}

\def\Bas{{\mathrm{Bas}}}
\def\End{{\mathrm{End}}}

\def\exter{{\textstyle\bigwedge}}

\def\exter{{\textstyle\bigwedge}}

{\theorembodyfont{\sl}
\newtheorem{theo}{Theorem}
\newtheorem{prop}[theo]{Proposition}
\newtheorem{cor}[theo]{Corollary}

\newtheorem{lem}[theo]{Lemma}}
{\theorembodyfont{\rm}

}

\begin{document}
\baselineskip=0.7cm

\vspace*{3cm}
\begin{center} 
{\large\bf ON THE NONCOMMUTATIVE GEOMETRY\\
\bigskip
OF THE ENDOMORPHISM ALGEBRA\\
\bigskip
OF A VECTOR BUNDLE} 
\end{center}
\vspace{1cm}

\begin{center}
Thierry MASSON\\
\vspace{0.3cm}
{\small Laboratoire de Physique Th\'eorique et Hautes
Energies\footnote{Laboratoire associ\'e au Centre National de la
Recherche Scientifique - URA D0063}\\
Universit\'e Paris XI, B\^atiment 211\\
91 405 Orsay Cedex, France\\
e-mail: masson@qcd.th.u-psud.fr}
\end{center}
\vspace{1cm}

\vspace{1cm}
\begin{abstract}
In this letter we investigate some aspects of the noncommutative differential 
geometry based on derivations of the algebra of endomorphisms of an oriented 
complex hermitian vector bundle. We relate it, in a natural way, to the geometry 
of the underlying principal bundle and compute the cohomology of its complex of  
noncommutative differential forms.
\end{abstract}

\vfill
\noindent L.P.T.H.E.-ORSAY 98/14\\
\vspace{0.5cm}

\newpage

\section{Introduction and notations}

In \cite{DVMA}, it was shown that the noncommutative geometry of the algebra of 
endomorphisms of an oriented complex hermitian vector bundle is very much like the 
ordinary geometry of the associated $SU(n)$-principal bundle. In particular, from 
the point of view of connections, this noncommutative algebra gives us interesting 
relations with the canonical Atiyah Lie algebroid associated to this oriented 
vector bundle. 

In this letter, we would like to proceed in the study of the noncommutative 
differential calculus of this endomorphisms algebra. In particular, we would like 
to make closer relations with the ordinary geometry of the principal bundle. Using 
ordinary technics in algebraic geometry, we will compute the cohomology of this 
noncommutative differential calculus.

The notations we will use in this letter are those of \cite{DVMA}, and we refer 
the reader to this paper for more details on the construction and properties of 
the objects we will introduce here.

Let $M$ be a regular finite dimensional smooth manifold. We denote by $E$ an 
oriented complex hermitian vector bundle of rank $n$ over $M$, by $P$ its 
$SU(n)$-principal frame bundle, and we introduce $\fA$ the algebra of 
endomorphisms of $E$, which is the set of sections of $\End(E)$. We denote by 
$(\Omega_\der(\fA), \hat{d})$ the noncommutative differential calculus based on 
derivations \cite{MDV1, DVMI, DVMA}.

\section{Some relations between $\Omega_\der(\fA)$ and $\Omega(P)$}

In this section we would like to give some structure properties of 
$\Omega_\der(\fA)$ which relates it to the ordinary differential calculus 
$\Omega(P)$ of $P$. 

Let us denote by $\cF(P)$ the (commutative) algebra of smooth functions on $P$ and 
by $A \mapsto A^v$ the map which sends any $A \in \tsu(n)$ into the associated 
vertical vector field on $P$. 

Let us introduce the algebra $\fB = \cF(P) \otimes M_n(\gC)$ of matrix valued 
functions on $P$.   Denote by $(\Omega_\der(\fB), \hat{d}) = (\Omega(P) \otimes 
\Omega_\der(M_n(\gC)), d + d')$ its differential calculus based on derivations 
\cite{DVKM2}. It is easy to see that $\cG = \{ A^v + ad_A\ / \ A \in \tsu(n) \}$ 
is a Lie subalgebra of $\der(\fB)$ isomorphic to $\tsu(n)$. This Lie subalgebra 
defines a Cartan operation of $\tsu(n)$ on $\Omega_\der(\fB)$, whose basic 
subalgebra we denote by $\Omega_{\der,\Bas} (\fB)$. Then one has:

\begin{prop}
$\Omega_\der(\fA) = \Omega_{\der,\Bas} (\fB)$
\end{prop}

\demo The proof is based on results on noncommutative quotient manifolds studied 
in \cite{MAS}. First, notice that $\fA$, as a set of section of an associated 
bundle of $P$, can be considered as the algebra $\{ b \in \fB\ / \ A^v a + [ A, a] 
= 0\ \forall \ A\in \tsu(n) \}$.  Now, define as in \cite{MAS}
\begin{eqnarray*}
\cZ_\der(\fA) &=& \{ \tX \in \der(\fB)\ / \ \tX \fA = 0 \} \\
\cN_\der(\fA) &=& \{ \tX \in \der(\fB)\ / \ \tX \fA \subset \fA \} 
\end{eqnarray*}
Then looking locally (in a trivialization of $P$) at the derivations of $\fB$, one 
sees that 
\begin{eqnarray*}
\cZ(\fA) &=& \fA \cap \cZ(\fB) \\
\der(\fA) &=& \cN_\der(\fA) / \cZ_\der(\fA) \\
\fA &=& \{ b \in \fB\ / \ \tX b = 0 \ \forall \ \tX\in \cZ_\der(\fA) \}
\end{eqnarray*}
where $\cZ(\fA)$ and $\cZ(\fB)$ are the center of the algebras $\fA$ and $\fB$ 
respectively. This makes $\fA$ into a noncommutative quotient manifold algebra of 
$\fB$ in the sense of \cite{MAS}. In order to prove the proposition, using 
Prop~V.1 in \cite{MAS} and the fact that $\underline{\Omega}_\der$ and 
$\Omega_\der$ coincide in this context \cite{DVMA}, it remains to show that the 
$\cZ(\fB)$-module induced by $\cN_\der(\fA)$ in $\der(\fB)$ is $\der(\fB)$ itself. 
As before, using local expressions of derivations, this can be checked 
easily.\findem

\medskip
As an example, let us consider a $SU(n)$-connection on $P$. Denote by $\omega$ its 
$1$-form on $P$. It was shown in \cite{DVMA} that there exists a corresponding 
noncommutative $1$-form $\alpha \in \Omega_\der(\fA)$. This form comes from a 
basic $1$-form in $\Omega_{\der,\Bas} (\fB)$, which is nothing but $\omega - i 
\theta$, where $\theta \in \Omega^1_\der(M_n(\gC))$ is the canonical $1$-form 
defined in \cite{DVKM1}. The basicity of this $1$-form is a consequence of 
properties of $\omega$ and $i\theta$, in particular the equivariance of $\omega$. 

\medskip
Now, notice that the commutative algebra $\cF(M)$ of smooth functions on $M$ and 
its de~Rham complex $\Omega(M)$ are the basic subalgebra of $\cF(P)$ and the basic 
subcomplex of $\Omega(P)$ for the operation of $\tsu(n)$ induced by $A \mapsto 
A^v$. This operation is itself the restriction of the operation of $\tsu(n)$ 
considered previously. Then, from this point of view, $\Omega_\der(\fA)$ is a 
natural generalization of $\Omega(M)$ containing informations on $P$.

\medskip
Moreover, this construction fits perfectly with the notion of noncommutative 
integration. It was shown in \cite{DVKM1} that such a notion exists on the 
noncommutative differential calculus $\Omega_\der(M_n(\gC))$. We denote by $\omega 
\in \Omega^{n^2 - 1}_\der(M_n(\gC)) \mapsto \int_\mathrm{nc} \omega \in \gC$ this 
noncommutative integral. This integral induces a map 
\begin{eqnarray*}
\Omega^{p, n^2 - 1}_\der(\fB) &\rightarrow & \Omega^p(P)\\
\omega &\mapsto & \int_\mathrm{nc} \omega 
\end{eqnarray*}
which satisfies:
\begin{prop}
1) If $\omega \in \Omega^{p, n^2 - 1}_\der(\fB)$ is horizontal (resp. invariant), 
then $\int_\mathrm{nc} \omega \in \Omega^p(P)$ is horizontal (resp. invariant) for 
the two operations defined above. 

2) Considering basic elements, this map defines a canonical noncommutative 
integration ``along the (noncommutative) fiber'' $\Omega_\der(\fA) \rightarrow 
\Omega(M)$.

3) This noncommutative integral is compatible with the differentials:
\[
\int_\mathrm{nc} \hat{d} \omega = d \int_\mathrm{nc} \omega
\]

4) This induces maps in cohomologies:
\begin{eqnarray*}
\int_{\mathrm{nc}\ \ast} : H^r( \Omega_\der(\fB) , \hat{d}) & \rightarrow & H^{r - 
(n^2 - 1)}(P) \\
\int_{\mathrm{nc}\ \ast} : H^r( \Omega_\der(\fA) , \hat{d}) & \rightarrow & H^{r - 
(n^2 - 1)}(M)
\end{eqnarray*}
\end{prop}

\demo 1) and 3) are straightforward computations using the precise definitions of 
the noncommutative integration over $\Omega_\der(M_n(\gC))$, the differentials, 
and the two operations.

2) and 4) are immediate consequences of 1) and 3).\findem

\medskip
The cohomology groups involved in 4) will be described in the next section where 
the computation of the cohomology of $\Omega_\der(\fA)$ is performed.

\medskip
The different relations between the various differential calculi can be summarized 
in the following diagram:
\[
\begin{diagram}
\node{\Omega_\der(\fA)} \arrow{s,l}{\int_\mathrm{nc}} 
\arrow{e,tb,J}{\mathrm{basic}}{\mathrm{elements}} \node{\Omega_\der(\fB)} 
\arrow{s,r}{\int_\mathrm{nc}} \\
\node{\Omega(M)} \arrow{e,tb,J}{\mathrm{basic}}{\mathrm{elements}}  
\node{\Omega(P)}
\end{diagram}
\]

\section{The cohomology of $\Omega_\der(\fA)$}

In the case when $\fA$ is a tensor product $\fA = \cF(M) \otimes M_n(\gC)$, the 
cohomology of $\Omega_\der(\fA)= \Omega(M) \otimes \Omega_\der(M_n(\gC))$ is known 
\cite{DVKM2}. It is just the tensor product of the cohomology of $\Omega(M)$ (the 
de~Rham cohomology) and the cohomology of $\Omega_\der(M_n(\gC))$ (we will detail 
this last cohomology in the following). 

In the general case, the cohomology of $\Omega_\der(\fA)$ can be computed using a 
slight modification of standard constructions in algebraic topology \cite{BT}. In 
this section, we perform this computation and show that the result is the same as 
in the tensor product case.

\medskip
Let $U$ be a open subset of $M$ such that the restriction of $\End(E)$ to $U$ is 
trivial. We make a choice of trivializations for any such open subset and we 
denote by $\fA(U)$ the trivialization of the restriction to $U$ of the algebra 
$\fA$. Then one has $\fA(U) \simeq C^\infty(U)\otimes M_n(\gC)$. Denote by $g_{UV} 
:  U \cap V \rightarrow SU(n)$ the transition functions. 

Consider now the presheaf $\cF$ defined by $U \mapsto \Omega_\der(\fA(U))$ where 
$U$ is any open subset of $M$ which trivializes $\End(E)$. For $V \subset U$, the 
inclusion map is defined to be
\begin{eqnarray*} 
i_U^V :  \Omega_\der(\fA(U)) &\rightarrow& \Omega_\der(\fA(V)) \\
\omega &\mapsto& \left( \omega_{\upharpoonright V} \right)^{g_{UV}} 
\end{eqnarray*}
where $\left( \omega_{\upharpoonright V} \right)^{g_{UV}}$ is the action of 
$g_{UV}$ (change of trivialization) on the restriction of $\omega$ to $V$. If 
$\omega = a_0 \hat{d} a_1 \dots \hat{d} a_p$, one has 
\[
\omega^{g_{UV}} = ( g_{UV}^{-1} a_0 g_{UV} ) \hat{d} ( g_{UV}^{-1} a_1 g_{UV} ) 
\dots \hat{d} ( g_{UV}^{-1} a_p g_{UV} )
\]
This action commute with $d$.

Now, let us take a good cover $\tU = \{ U_\alpha\}_{\alpha \in I}$ of $M$ indexed 
by an ordered set $I$ and such that over any $U_\alpha$, $\End(E)$ is trivialized. 
For convenience, on any $U_{\alpha_1 \dots \alpha_q} = U_{\alpha_1} \cap \dots 
\cap U_{\alpha_q}$ the trivialization is chosen to be the restriction of the 
trivialization of $U_{\alpha_q}$.

We now define a noncommutative version of the double \v{C}ech-de~Rham complex 
associated with this presheaf. For $p \geq 0$ and $q \geq 0$, consider the vector 
spaces 
\[ C^{p,q}( \tU, \cF) = \prod_{\alpha_0 < \dots < \alpha_p} \Omega^q_\der( \fA( 
U_{\alpha_0 \dots \alpha_p})) \]
An element $\omega \in C^{p,q}( \tU, \cF)$ is a collection of noncommutative 
$q$-form in $\omega_{\alpha_0 \dots \alpha_p} \in \Omega^q_\der( \fA( U_{\alpha_0 
\dots \alpha_p}))$. 

Define the differential $\delta :  C^{p, q}( \tU, \cF) \rightarrow  C^{p+1, q}( 
\tU, \cF)$ by
\[
(\delta \omega)_{\alpha_0  \dots  \alpha_{p+1}} = \sum_{i = 0}^p (-1)^i 
\omega_{\alpha_0  \dots \alpha_{i-1} \alpha_{i+1} \dots \alpha_{p+1}} + (-1)^{p+1} 
\omega_{\alpha_0  \dots  \alpha_p}^{g_{\alpha_p \alpha_{p+1}}} 
\]
Using the properties of the transition functions it is easy to verify that 
$\delta^2 = 0$. The noncommutative differential $\hat{d}$ is of degree $(0, 1)$ on 
this double complex and satisfies $\hat{d} \delta = \delta \hat{d}$. On the total 
complex of this double complex, we introduce the differential $D = \delta + (-1)^p 
\hat{d}$, with $D^2 = 0$. 

For $p=-1$, define $C^{-1, q}( \tU, \cF)$ to be $\Omega_\der^q(\fA)$, and $\delta 
: C^{-1, q}( \tU, \cF) \rightarrow C^{0, q}( \tU, \cF)$ to be the restriction map. 
Then $\delta^2 = 0$ also holds on $C^{-1, q}( \tU, \cF)$.

\begin{lem}
The following sequence is exact :
\[
\dgARROWLENGTH=1.5em
\begin{diagram}
\node{0} \arrow{e} \node{ C^{-1,\ast}( \tU, \cF)} \arrow{e,t}{\delta} \node{ 
C^{0,\ast}( \tU, \cF)} \arrow{e,t}{\delta} \node{C^{1,\ast}( \tU, \cF)} 
\arrow{e,t}{\delta} \node{\dots} 
\end{diagram}
\]
\end{lem}

\demo The exactitude at $C^{-1,\ast}( \tU, \cF)$ is trivial. Because $\delta^2 = 
0$, one has only to show that if $\omega \in C^{p,\ast}( \tU, \cF)$, with $\delta 
\omega = 0$, then there exists $\eta \in C^{p-1,\ast}( \tU, \cF)$ such that 
$\delta \eta = \omega$. Introduce $\{ \rho_\alpha \}$ a partition of unity 
subordinate to the good cover $\tU$. For $\alpha_0 < \dots < \alpha_{p-1}$, define 
\[
\eta_{\alpha_0 \dots \alpha_{p-1}} = \sum_{{\alpha < \alpha_{p-1}} \atop {\alpha 
\neq \alpha_0, \dots, \alpha_{p-2}} } \rho_\alpha \omega_{\alpha \alpha_0 \dots 
\alpha_{p-1}} + \sum_{ \alpha > \alpha_{p-1} } \rho_\alpha \omega_{\alpha \alpha_0 
\dots \alpha_{p-1}}^{ g_{\alpha \alpha_{p-1} }}
\]
where, to simplify this expression, we make use of the notation $\omega_{\dots 
{\alpha_i} \dots {\alpha_j} \dots} = - \omega_{\dots {\alpha_j} \dots {\alpha_i} 
\dots}$. Note that $\rho_\alpha \omega_{\alpha \alpha_0 \dots \alpha_{p-1}} \in 
\Omega_\der(\fA(U_{\alpha \alpha_0 \dots \alpha_{p-1}}))$. A straightforward 
computation shows then that with $\delta\omega = 0$, one has $\delta \eta = 
\omega$.\findem

\medskip
Using general arguments on double complexes, this lemma gives us a noncommutative 
version of the Generalized Mayer-Vietoris Principle: 

\begin{cor}
The cohomology of $(\Omega_\der(\fA), \hat{d})$ is the cohomology of the total 
complex $(C(\tU, \cF), D)$.
\end{cor} 

\medskip
Consider now the spectral sequence $\{E_r, d_r\}$ associated with the double 
complex $C(\tU, \cF)$, induced by the filtration 
\[
F^p C(\tU, \cF) = \oplus_{s \geq p} \oplus_{q \geq 0} C^{s,q}(\tU, \cF)
\]
Then by standard argument, one knows that
the first term of this spectral sequence is
\[
E_1^{p,q} = H_{\hat{d}}^{p,q} = C^p(\tU, \cH^q)
\] 
where $\cH^q$ is the presheaf $\cH^q(U) = H^q(\Omega_\der(\fA(U)), \hat{d})$ and 
that the second term is
\[
E_2^{p,q} = H_{\delta}^p( \tU, \cH^q)
\]
This spectral sequence converges to $H(C(\tU, \cF), D) = H(\Omega_\der(\fA), 
\hat{d})$.

The cohomology groups $H^\ast(\Omega_\der(\fA(U)), \hat{d})$ have been computed 
\cite{MDV1, DVKM1, DVKM2}. When $U$ is diffeomorphic to $\gR^m$, one has
\[
H^\ast(\Omega_\der(\fA(U)), \hat{d}) = H^\ast( \Omega_\der(M_n(\gC)), d')
\]
and $H^\ast( \Omega_\der(M_n(\gC)), d')$ is isomorphic to the cohomology of the 
Lie algebra $\tsl(n,\gC)$. Then $H^\ast( \Omega_\der(M_n(\gC)), d') = \left( 
\exter \tsl(n,\gC)^\ast \right)_{\mathrm{Inv}}$, where Inv denotes the invariant 
elements for the canonical Lie derivation. This algebra can also be considered as 
the graded commutative algebra freely generated by $c^n_{2 r - 1}$ in degree 
$2r-1$ with $r \in \{ 2,3,\dots, n\}$ and where the $c^n_{2 r - 1}$ are the 
primitive elements of $\left( \exter \tsl(n,\gC)^\ast \right)_{\mathrm{Inv}}$.
 
Any element in $\left( \exter \tsl(n,\gC)^\ast \right)_{\mathrm{Inv}}$ is then 
invariant by the action of $SU(n)$. So, the inclusion map is the identity and then 
the presheaf $\cH^q$ is a constant presheaf. The cohomology of $E_1$ reduces to 
the \v Cech cohomology of this constant presheaf. Now, using again the structure 
of $H^\ast( \Omega_\der(M_n(\gC)), d')$, it is easy to see that the differentials 
$d_2, \dots, d_r$ are all zero. We have then proven:

\begin{prop}
$E_\infty = E_2 = H^\ast(M) \otimes \left( \exter \tsl(n,\gC)^\ast 
\right)_{\mathrm{Inv}}$
\end{prop}

This proposition tells us that the cohomology of $\Omega_\der(\fA)$ is far more 
simpler that the cohomology of the underlying principal bundle. 

Notice now that there is a canonical inclusion of algebras
\[
\left( \exter \tsl(n,\gC)^\ast \right)_{\mathrm{Inv}} \hookrightarrow \left( 
\exter \tsl(n+1,\gC)^\ast \right)_{\mathrm{Inv}}
\]
given by $c^n_{2 r - 1} \mapsto c^{n+1}_{2 r - 1}$ for $r \in \{ 2,3,\dots, n\}$. 
Taking the inductive limit, one defines $\left( \exter \tsl(\infty,\gC)^\ast 
\right)_{\mathrm{Inv}}$, which is the graded commutative algebra freely generated 
by  elements $c_{2r-1}$ in degree $2r-1$. Then, the cohomology of 
$\Omega_\der(\fA)$ can be embedded in $H(M) \otimes \left( \exter 
\tsl(\infty,\gC)^\ast \right)_{\mathrm{Inv}}$, the right hand side being now 
completely independent of the oriented hermitian vector bundle $E$.

\section{Conclusions}

In this letter, we have added a few more arguments to those given in \cite{DVMA}, 
to propose the algebra $\fA$ equipped with its noncommutative differential 
calculus based on derivations as a possible replacement of the principal bundle 
$P$. Indeed, we have shown that this noncommutative geometry of $\fA$ is strongly 
and very naturally related to the ordinary geometry of $P$. Then $\fA$ can be used 
in place of $P$, if one replaces the differential calculus $\Omega(P)$ by 
$\Omega_\der(\fA)$, the connection $1$-form  $\omega$ on $P$ by the associated 
noncommutative $1$-form $\alpha$ introduced in \cite{DVMA} (which is only 
subjected to a ``vertical'' condition), the notion of associated vector bundle by 
the notion of (left-)module over $\fA$.

From a physical point of view, this noncommutative geometry is more interesting 
because, as was pointed out in \cite{DVKM2, DVMA}, it contains not only ordinary 
Yang-Mills fields, but also other fields which look very much like Higgs fields. 
On the other hand, from a mathematical point of view, this geometry looks simpler, 
in particular in cohomology.

\section*{Acknowledgments}

 We would like to thank Michel Dubois-Violette for very helpful discussions and 
his kind interest.

\end{document}